\documentclass[a4paper]{cas-sc}
\usepackage[numbers]{natbib}
\usepackage{url}
\usepackage{doi}
\usepackage{amsmath, amssymb, amsthm}
\usepackage{graphicx}
\usepackage{booktabs}
\usepackage[ruled, linesnumbered]{algorithm2e}
\usepackage{hyperref}
\hypersetup{colorlinks=true, linkcolor=blue, citecolor=blue}
\usepackage{caption}
\providecommand{\keywords}[1]{}

\usepackage{fancyhdr}
\fancypagestyle{first}{
  \fancyhf{}
  
  \fancyfoot[L]{}
  \fancyfoot[R]{\thepage}
}

\begin{document}

\title{\fontsize{20}{24}\selectfont GAC-PINN: Geometry-Adaptive and Constraint-Enhanced Physics-Informed Neural Networks}

\author[1]{Yanxin Zhang}
\author[1]{Yong Zhang}
\author[1]{Houbiao Li\corref{cor1}}[orcid=0000-0002-7268-307X]
\ead{lihoubiao0189@163.com}

\affiliation[1]{
    organization={University of Electronic Science and Technology of China},
    city={Chengdu},
    postcode={611731},
    state={Sichuan},
    country={China}
}

\cortext[cor1]{Corresponding author.}

\begin{abstract}
For systems with steep gradients, sharp interfaces, or severe spatio-temporal coupling, Physics-informed neural networks (PINNs) suffer from spectral bias, geometric inflexibility, and boundary constraint conflicts, which undermine accuracy and convergence. To overcome these issues, we propose a geometry-adaptive and constraint-enhanced PINN (GAC-PINN). The framework comprises four components: a gradient-driven adaptive grid mapping (AGM) for diffeomorphic point concentration with Jacobian regularization, an adaptive bandwidth hard-constraint ansatz with spatially-varying boundary transition widths, a Gaussian Fourier feature mapping as a spectral preconditioner to further enhance high-wavenumber representation, and an operator-aware router that automatically selects the appropriate hard-constraint construction based on whether the governing PDE contains temporal derivatives. An AGM callback mechanism and a three-stage training strategy ensure stable coordination. Benchmarks including the viscous Burgers equation, a sharp-peaked 2D Poisson problem, and the Allen-Cahn phase-transition equation show that GAC-PINN attains relative \(L^2\) errors of \((1.747\pm0.450)\times10^{-4}\), \((2.868\pm0.947)\times10^{-5}\), and \((1.756\pm0.712)\times10^{-3}\), respectively, consistently outperforming the baselines. Ablation studies further reveal that AGM alone yields a substantially lower error than residual-based adaptive refinement (RAR), while RAR becomes beneficial only when combined with FFM, demonstrating a context-dependent module interaction. Convergence analysis verifies rapid error reduction and saturation with increasing resolution, establishing a practical adaptive framework for high-fidelity simulation of problems with localized sharp features in applied mechanics and computational physics.
\end{abstract}


\shorttitle{Geometry-Adaptive and Constraint-Enhanced PINNs}
\shortauthors{Zhang et al.}

\maketitle


\section{Introduction}

Accurate simulation of physical fields with steep gradients, sharp interfaces, and strong nonlinearities is critically important in numerous engineering disciplines. In aerospace engineering, the precise capture of shock waves and boundary layers is critical for predicting aerodynamic drag, heat flux, and structural integrity of high-speed vehicles; in electronic packaging, reliable thermal management requires high-fidelity solutions of the heat equation in the presence of localized hot spots and material interfaces; and in materials science, the modeling of phase separation and grain growth during alloy solidification hinges on the accurate resolution of propagating phase interfaces under extreme curvature-driven dynamics. These diverse applications share a common mathematical challenge: the efficient and accurate numerical solution of partial differential equations (PDEs) whose solutions exhibit localized, multi-scale features that impose prohibitive resolution requirements on conventional mesh-based methods.

Physics-informed neural networks (PINNs), introduced by Raissi et al. \cite{raissi2019physicsinformed}, have emerged as a transformative paradigm for PDE solving. By embedding physical laws directly into the training loss and leveraging automatic differentiation (AD) to compute PDE residuals, PINNs provide a mesh-free framework that circumvents the discretization and meshing burdens of traditional numerical methods. This elegant formulation has enabled promising results in various forward and inverse problems.

However, when deployed on the aforementioned problems with strong localized gradients and multiscale dynamics, standard PINNs suffer from three interrelated bottlenecks that severely compromise accuracy and convergence robustness:

\begin{itemize}
\item \textbf{Spectral bias:} Deep fully-connected networks prefer to learn low-frequency components due to the eigenvalue decay of the neural tangent kernel (NTK), causing pronounced oscillations near shocks or boundary layers.

\item \textbf{Geometric inflexibility:} Fixed uniform collocation points fail to dynamically adapt to evolving solution features, leading to redundant sampling in smooth regions and insufficient resolution in high-gradient zones.

\item \textbf{Boundary/initial condition conflicts:} Penalty-based soft constraints induce detrimental gradient competition between PDE residuals and boundary losses, while conventional fixed hard-constraint ansatz may introduce derivative contamination in high-order AD computations.

\end{itemize}

Existing works have achieved significant progress by tackling individual bottlenecks. 
For instance, Fourier feature mappings \cite{tancik2020fourier} mitigate spectral bias, 
residual-based adaptive sampling \cite{wu2023comprehensive} optimizes collocation point 
distribution, and coordinate transformation techniques \cite{gao2021phygeonet} enhance 
geometric flexibility. However, these strategies are often applied as independent 
modules or in a sequential pipeline, leaving room for a more synergistic integration 
that jointly considers geometric adaptation and spectral preconditioning within a 
unified training framework. In this work, we propose such an integration, with a 
particular focus on a gradient-driven diffeomorphic mapping (AGM) that dynamically 
adapts collocation points based on real-time physical gradients. This mapping is 
designed to be compatible with and complementary to existing spectral embedding 
techniques, thereby providing a cohesive framework rather than a simple aggregation 
of disjoint components.

Building on this integrated perspective, this work proposes a Geometry-Adaptive and Constraint-Enhanced PINN (GAC-PINN). The framework consists of four components integrated within a unified training pipeline that directly target the three bottlenecks: an operator-aware router that selects the appropriate hard-constraint construction based on the presence of temporal derivatives; a gradient-driven adaptive grid mapping (AGM) that diffeomorphically concentrates collocation points in high-residual regions while preventing geometric degeneration; a Gaussian Fourier feature mapping (FFM) that reshapes the NTK spectrum to mitigate low-frequency bias; and an adaptive hard-constraint ansatz with a learnable bandwidth and a stop-gradient operator that facilitates stable high-order automatic differentiation. A dedicated three-stage training strategy ensures stable coordination among these components. The main contributions of this work are threefold:

\begin{itemize}

\item \textbf{Methodologically,} we establish a new PINN framework that achieves tight integration of geometric adaptation, spectral preconditioning and constraint enforcement within a single optimization pipeline;

\item \textbf{Algorithmically,} we propose a closed-loop manifold evolution mechanism——AGM driven by real-time physical gradients, distinct from static coordinate transformations or discrete resampling;

\item \textbf{Experimentally,} extensive benchmarks on the viscous Burgers equation, a sharp 2D Poisson problem, the Allen-Cahn phase-field equation, and the 2D Navier-Stokes equations demonstrate that GAC-PINN achieves competitive and improved accuracy compared to re-implemented baselines under aligned settings, while ablation and computational analyses confirm the synergistic contributions of each module and competitive training efficiency.
\end{itemize}

The remainder of this paper is organized as follows. Section 2 reviews the relevant preliminaries, including PINNs, NTK theory, and diffeomorphic mappings. Section 3 presents the detailed architecture and algorithmic implementation of GAC-PINN. Section 4 describes the experimental setup. Section 5 reports comprehensive numerical results and comparisons. Finally, Section 6 provides discussion on the implications and limitations of the proposed approach, also outlines future research directions.

\section{Methodology}

\subsection{Physics-informed neural networks}

Consider a spatiotemporal domain \(\Omega \times [0, T]\), with spatial coordinates \(x \in \mathbb{R}^d\) and time \(t \in [0, T]\). Let \(X = (x, t)\) denote the spatiotemporal coordinate. The general PDE initial-boundary value problem is expressed as:
\begin{align}
\mathcal{P}[u](X) &= 0, \quad X \in \Omega \times (0, T], \\
\mathcal{B}[u](X) &= 0, \quad X \in \partial\Omega \times (0, T], \\
\mathcal{I}[u](X) &= 0, \quad x \in \Omega, \ t = 0,
\end{align}
where \(\mathcal{P}\) denotes the differential operator, \(\mathcal{B}\) and \(\mathcal{I}\) represent boundary and initial operators, respectively, and \(u(X)\) is the unknown physical field.A standard PINN \cite{raissi2019physicsinformed} approximates the unknown solution by a fully connected neural network \(u(X; \Theta)\) with trainable parameters \(\Theta\). The total loss function is a weighted sum of three components:
\begin{equation}
\mathcal{L}_{PINN} = \lambda_{pde}\mathcal{L}_{pde} + \lambda_{bc}\mathcal{L}_{bc} + \lambda_{ic}\mathcal{L}_{ic},
\end{equation}
where the PDE residual loss quantifies the deviation of the governing equation at collocation points:
\begin{equation}
\mathcal{L}_{pde} = \frac{1}{N_{pde}}\sum_{i=1}^{N_{pde}} \left\| \mathcal{P}[u](X_i) \right\|_2^2,
\end{equation}
and the boundary and initial losses are defined analogously. The weighting coefficients \(\lambda_{pde}\), \(\lambda_{bc}\), and \(\lambda_{ic}\) balance the competing loss terms, and the network parameters are optimized via back-propagation to minimize this composite objective.

\subsection{Neural tangent kernel theory}

The neural tangent kernel (NTK) provides a theoretical tool for analyzing the training dynamics of infinite-width neural networks. For a fully connected network \(u(X; \Theta)\), the NTK is defined as \cite{jacot2018neural}:
\begin{equation}
\Theta_{NTK}(X, X') = \left\langle \frac{\partial u(X; \Theta)}{\partial \Theta}, \frac{\partial u(X'; \Theta)}{\partial \Theta} \right\rangle.
\end{equation}

Under gradient descent training with an infinitesimal learning rate, the NTK remains approximately constant during training, i.e., \(\Theta_t \approx \Theta_0\). For a training set of \(N\) samples,the empirical NTK matrix \(\Theta \in \mathbb{R}^{N \times N}\) is symmetric positive semi-definite and admits the eigendecomposition \cite{wang2022pinns}:
\begin{equation}
\Theta = \sum_{k=1}^{N} \lambda_k \phi_k \phi_k^T,
\end{equation}
with eigenvalues $\lambda_1 \ge \lambda_2 \ge \cdots \ge \lambda_N \ge 0$ and corresponding orthogonal eigenvectors $\phi_k$. Projecting the residual $e(t) = u(t) - y$ onto this eigenbasis yields modal coefficients $\hat{\epsilon}_k(t) = \phi_k^\top e(t)$. Due to eigenvector orthogonality, the high-dimensional training dynamics decouple into $N$ independent scalar ordinary differential equations:
\begin{equation}
\frac{d\hat{\epsilon}_k(t)}{dt} = -\eta \lambda_k \hat{\epsilon}_k(t), \label{eq:ntk_ode}
\end{equation}
which admits the closed-form solution
\begin{equation}
\hat{\epsilon}_k(t) = \hat{\epsilon}_k(0) e^{-\eta \lambda_k t}. \label{eq:ntk_solution}
\end{equation}

This theory reveals that high-frequency modes corresponding to small eigenvalues converge at exponentially slower rates than low-frequency modes, mathematically establishing the spectral bias of standard PINNs — a limitation that motivates the Fourier feature embedding introduced in Section~\ref{subsec:ffm}.

\subsection{Diffeomorphic coordinate transformation}
\label{subsec:diffeo}

\subsubsection{Definition of diffeomorphic mapping}
\label{subsubsec:diffeo_def}

A mapping $\mathcal{M}: \Omega_\Xi \to \Omega_X$ from the computational coordinates $\Xi=(\xi,\eta)$ to the physical coordinates $X=(x,t)$ is \emph{diffeomorphic} \cite{gao2021phygeonet} if it satisfies: (i) \textbf{bijectivity} — each computational point maps to a unique physical point and vice versa; (ii) \textbf{infinite differentiability} — all derivatives exist and are continuous throughout the domain; and (iii) \textbf{invertible Jacobian} — the Jacobian matrix $\mathbf{J} = \partial X / \partial \Xi$ satisfies $\det(\mathbf{J}) > 0$ everywhere, preventing coordinate folding or local volume collapse.

\subsubsection{Volume transformation and sampling density}
\label{subsubsec:volume}

Unlike conventional fixed-grid methods where collocation points are prescribed directly in the physical domain, this approach adopts a computational-coordinate perspective commonly used in r-adaptivity. Let $\Xi \in \Omega_\Xi$ denote the computational coordinates with a uniform reference distribution, and let $\mathcal{M}: \Omega_\Xi \to \Omega_X$ be a diffeomorphic mapping from the computational domain to the physical domain:
\begin{equation}
X = \mathcal{M}(\Xi), \quad \Xi \in \Omega_\Xi \subset \mathbb{R}^d, \quad X \in \Omega_X \subset \mathbb{R}^d. 
\end{equation}

The Jacobian matrix of this mapping is
\begin{equation}
\mathbf{J}(\Xi) = \frac{\partial X}{\partial \Xi} = 
\begin{bmatrix}
\dfrac{\partial X_1}{\partial \Xi_1} & \dfrac{\partial X_1}{\partial \Xi_2} & \cdots & \dfrac{\partial X_1}{\partial \Xi_d} \\[10pt]
\dfrac{\partial X_2}{\partial \Xi_1} & \dfrac{\partial X_2}{\partial \Xi_2} & \cdots & \dfrac{\partial X_2}{\partial \Xi_d} \\[10pt]
\vdots & \vdots & \ddots & \vdots \\[10pt]
\dfrac{\partial X_d}{\partial \Xi_1} & \dfrac{\partial X_d}{\partial \Xi_2} & \cdots & \dfrac{\partial X_d}{\partial \Xi_d}
\end{bmatrix}, \quad 
\det(\mathbf{J}(\Xi)) > 0, \quad \forall \Xi \in \Omega_\Xi,
\end{equation}
where the positivity condition ensures bijectivity and prevents local collapse.

By the change-of-variables formula, the physical and computational volume elements satisfy
\begin{equation}
dV_X = |\det(\mathbf{J}(\Xi))| \, dV_\Xi. 
\end{equation}

Let $\rho_\Xi(\Xi) \equiv \text{const}$ denote the uniform sampling density in the computational domain. Since the total number of collocation points is invariant under $\mathcal{M}$, the physical-space density $\rho_X(X)$ satisfies
\begin{equation}
\rho_X(X) \, dV_X = \rho_\Xi(\Xi) \, dV_\Xi, 
\end{equation}
which yields
\begin{equation}\label{eq:density_mapping}
\rho_X(X) = \frac{\rho_\Xi(\Xi)}{|\det(\mathbf{J}(\Xi))|}. 
\end{equation}

Formula (\ref{eq:density_mapping}) reveals the core principle of adaptive node refinement: the physical sampling density is inversely proportional to the local Jacobian determinant. Regions with large $\det(\mathbf{J})$ thus receive denser collocation points, concentrating computational resources in steep-gradient areas, while regions with small $\det(\mathbf{J})$ are automatically coarsened. In essence, the Jacobian determinant acts as a magnification factor that stretches the uniform computational grid into a solution-adaptive non-uniform physical grid.

Motivated by the classical equidistribution principle \cite{deboor1973good,budd2009adaptivity}, we seek a mapping that minimizes the coefficient of variation of the weighted volumes $\omega(X)|\det(\mathbf{J}(\Xi))|$ across all collocation points, where $\omega(X)$ is a monitor function reflecting the local gradient magnitude. The AGM module introduced in Section\ref{subsec:adc} explicitly performs this minimization, ensuring that high-gradient regions receive a proportionally larger share of computational resources.

\section{The GAC-PINN framework}

This work aims to develop a high-precision, adaptive numerical solver for PDEs with steep gradients, strong nonlinearities, and multiscale features. This methodology makes original contributions to \textit{adaptive discretization} and integrates it with \textit{spectral preconditioning} into a unified PINN framework as shown in Fig.~\ref{fig:GAC-PINN}. Among these, the hard-constraint ansatz is employed as an engineering enhancement to stabilize high-order differentiation.

\begingroup
\centering
\includegraphics[width=0.8\textwidth]{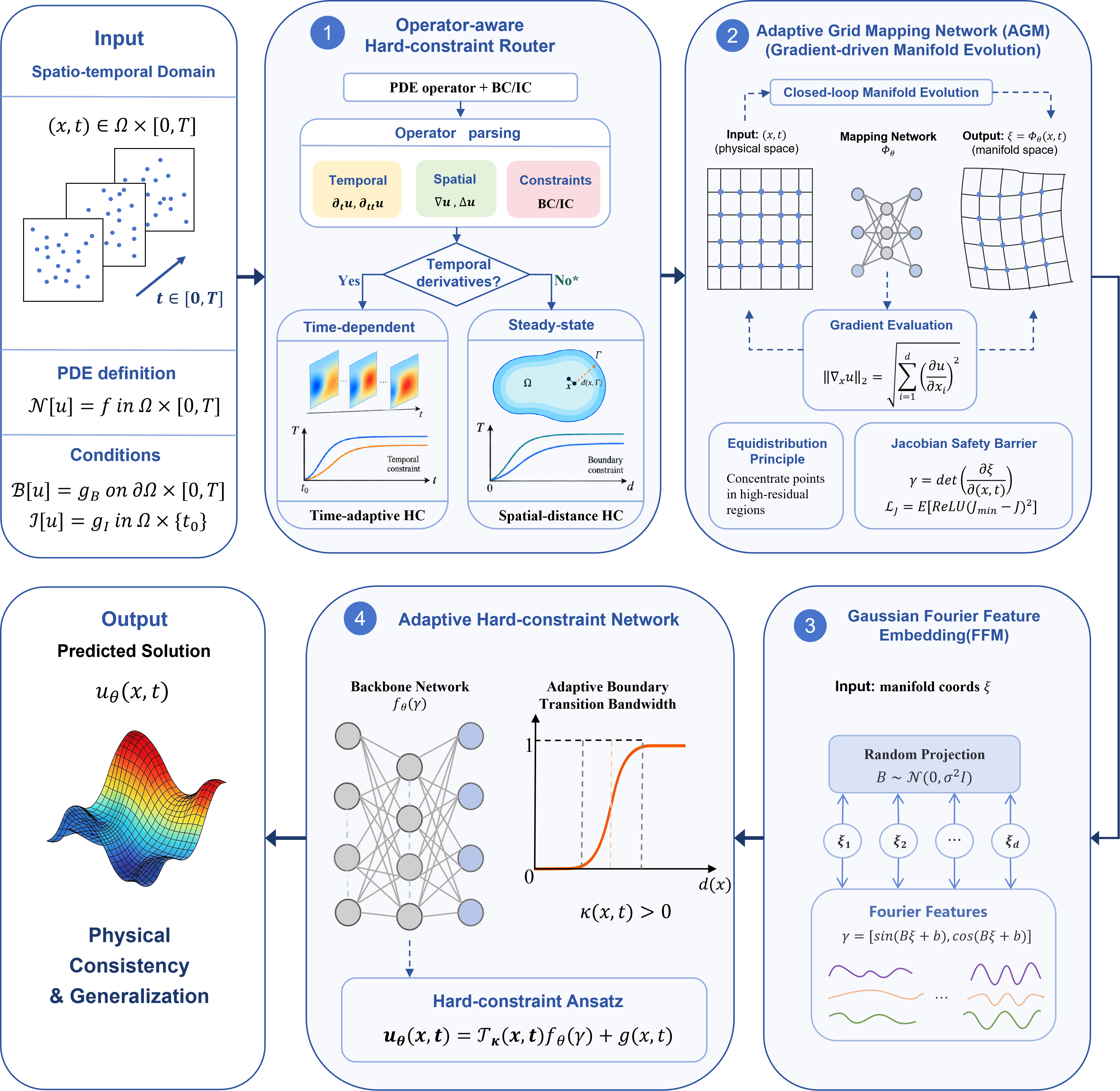}
\captionof{figure}{GAC-PINN Framework}
\label{fig:GAC-PINN}
\endgroup
\medskip

\subsection{Overall architecture}

Within this architecture, the neural network serves as the flexible functional approximator, while four components collectively address the above bottlenecks:

\begin{itemize}
\item \textbf{Operator-aware hard-constraint router:} Classifies the governing PDE by its dependence on temporal derivatives and selects the appropriate hard-constraint construction (time-adaptive bandwidth for evolution equations, pure spatial distance for steady-state problems).
\item \textbf{Adaptive grid mapping network (AGM):} Driven by physical gradients through the equidistribution principle, dynamically concentrates collocation points toward high-gradient regions with a Jacobian safety barrier against degeneration.
\item \textbf{Gaussian Fourier feature mapping (FFM):} Reshapes the NTK spectral distribution through Gaussian random projection, breaking the low-frequency bias.
\item \textbf{Adaptive hard-constraint ansatz:} Employs a learnable boundary transition bandwidth to facilitate stable high-order differentiation.
\end{itemize}

\subsection{Operator-aware hard-constraint routing}
\label{subsec:router}

The construction of the hard-constraint ansatz depends on the temporal characteristics of the governing PDE. Evolution equations benefit from a time-adaptive bandwidth that gradually releases the initial condition as $t$ increases, whereas steady-state problems require only a spatial distance function. To accommodate both classes within a unified framework without manual reconfiguration, we introduce an operator-aware router that classifies the governing operator and selects the corresponding hard-constraint construction.

The classification is based on the operator's dependence on temporal derivatives:
\begin{equation}
\chi_{\mathcal{P}} =
\begin{cases}
0, & \text{if } \partial\mathcal{P}/\partial(\partial_t u) \equiv 0 \quad \text{(steady-state)},\\[4pt]
1, & \text{if } \partial\mathcal{P}/\partial(\partial_t u) \not\equiv 0 \quad \text{(time-dependent)}.
\end{cases}
\label{eq:router}
\end{equation}

The discriminant $\chi_{\mathcal{P}}$ is determined symbolically from the PDE definition and incurs no additional computational cost. Based on $\chi_{\mathcal{P}}$, the framework selects the corresponding hard-constraint ansatz described in Section~\ref{subsec:ahc}:

\begin{itemize}
\item If $\chi_{\mathcal{P}} = 1$ (time-dependent), the \emph{time-adaptive bandwidth ansatz} is applied. A learnable gating network predicts a spatially-varying bandwidth $\kappa_{adaptive}(X)$ from the input coordinates, and the hard constraint gradually releases the initial condition as $t$ increases.
\item If $\chi_{\mathcal{P}} = 0$ (steady-state), a \emph{pure spatial distance function} is used, and all temporal dependence is removed.
\end{itemize}

This routing mechanism ensures that the same framework handles both elliptic and evolution problems without manual reconfiguration, while preserving the exact satisfaction of boundary and initial conditions in both cases.

\subsection{Gradient-driven adaptive grid mapping network (AGM)}
\label{subsec:agm}

\subsubsection{Gradient-driven feedback monitoring mechanism}

The driving force for manifold adaptation originates from the real-time first-order spatial derivatives of the physical network output $u_{NN}$ with respect to spatial coordinates $x$. At the current training iteration, automatic differentiation extracts the spatial gradient magnitude:
\begin{equation}
\|\nabla_x u_{NN}\|_2 = \sqrt{\sum_{i=1}^d \left( \frac{\partial u_{NN}}{\partial x_i} \right)^2}.
\end{equation}

To construct a dimensionless mapping from physical gradients to spatial grid density, we define a dynamic feedback indicator \(\omega(X)\) with global adaptive scaling:
\begin{equation}
\omega(X) = 1.0 + \alpha_{grad} \times \frac{\|\nabla_x u_{NN}\|_2}{\max_{\Omega_x} \|\nabla_x u_{NN}\|_2 + 10^{-8}},
\end{equation}
with gradient amplification coefficient \(\alpha_{grad} = 4.0\). This indicator normalizes by the current maximum gradient magnitude, strictly confining \(\omega(X)\) to the interval \([1.0, 5.0]\).

\subsection{Adaptive diffeomorphic coordinate transformation}
\label{subsec:adc}

Given that primary high-frequency physical features evolve strongly along spatial directions, this framework maintains absolute independence of the time coordinate to prevent spatial distortion from interfering with the temporal evolution. The spatial nonlinear transformation takes the form:
\begin{equation}
(x,y) = (\xi,\eta) + \alpha \cdot \tanh\left(\mathcal{A}_{GM}(\xi,\eta;\Theta_{MAP})\right), 
\end{equation}
where $\mathcal{A}_{GM}$ is a multilayer perceptron with trainable parameters $\Theta_{MAP}$, and $\alpha = 2.0$ controls the maximum local stretching magnitude.

\begingroup
\centering
\includegraphics[width=0.8\textwidth]{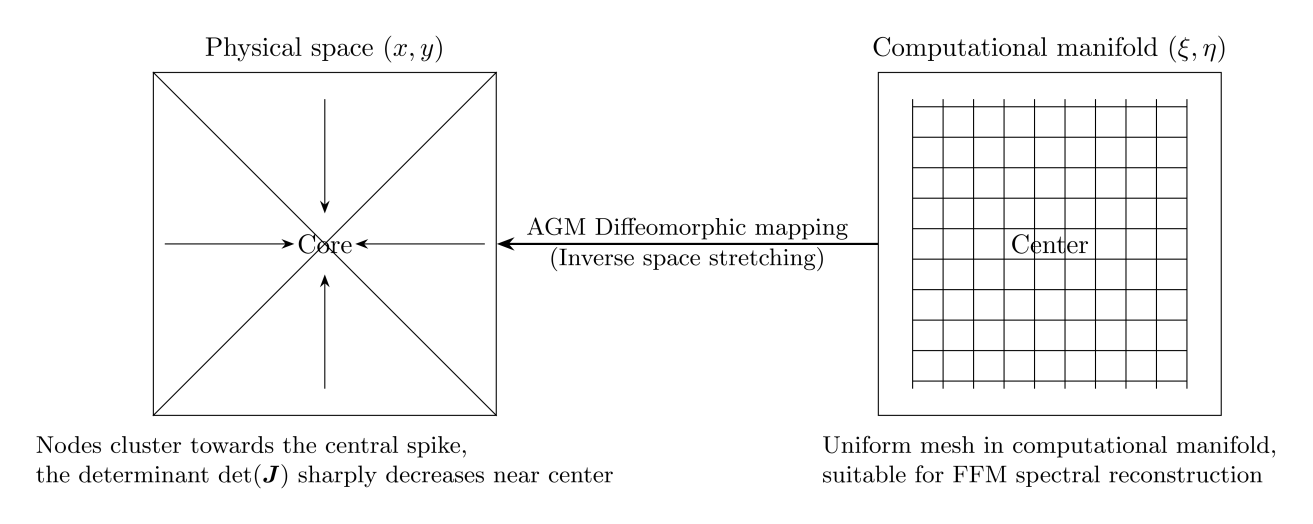}
\captionof{figure}{AGM geometric coordinate transformation}
\label{fig:agm}
\endgroup
\medskip

As illustrated in Fig.~\ref{fig:agm}, this transformation establishes a diffeomorphic mapping from the uniform computational manifold $(\xi, \eta)$ to the physical space $(x, y)$. It acts as a space-stretching operator: the uniform lattice in the computational domain is projected onto the physical domain, yielding adaptive node clustering toward high-gradient regions (e.g., the central spike). In these regions, the Jacobian determinant $\det(\mathbf{J})$ increases, reflecting the local volumetric magnification that concentrates collocation points. Conversely, smooth regions correspond to smaller $\det(\mathbf{J})$ and are automatically coarsened. Meanwhile, the computational manifold retains a perfectly uniform mesh, which is an essential prerequisite for subsequent FFM spectral reconstruction to guarantee orthogonality and numerical conditioning of the basis functions.

To quantitatively enforce this adaptive redistribution, the Jacobian matrix $\mathbf{J}$ (from computational to physical coordinates) is defined as:
\begin{equation}
\mathbf{J} = \frac{\partial(x, y)}{\partial(\xi, \eta)} =
\begin{bmatrix}
\dfrac{\partial x}{\partial \xi} & \dfrac{\partial x}{\partial \eta} \\[8pt]
\dfrac{\partial y}{\partial \xi} & \dfrac{\partial y}{\partial \eta}
\end{bmatrix}.
\end{equation}

Motivated by the equidistribution principle, we construct a weighted volumetric measure $V_i = \omega(\mathbf{X}_i) \cdot \det(\mathbf{J}_i)$ to characterize the local solution complexity. The AGM parameters are then optimized by a composite loss function that jointly minimizes the coefficient of variation of these volumes and enforces a Jacobian safety barrier:
\begin{align}
\mathcal{L}_{AGM} &= 2.0 \cdot \mathcal{L}_{equi} + \lambda_{barrier} \cdot \mathcal{L}_{guard},  \\
\end{align}
where
\begin{align}
\mathcal{L}_{equi} &= \frac{\sqrt{ \frac{1}{N} \sum_{i=1}^N (V_i - \bar{V})^2 }}{\bar{V} + \epsilon_0},  
\end{align}
Minimizing $\mathcal{L}_{equi}$ drives the physical grid volumes to conform to the equidistribution criterion, thereby automatically inducing local refinement in high-gradient regions, where the  computational lattice and $\det(\mathbf{J})$ increases to concentrate collocation points. Conversely, in smooth regions  where $\det(\mathbf{J})$ is small, the grid is automatically coarsened. The safety barrier, with a lower bound $J_{safe\_min} = 0.02$, imposes a quadratic penalty when the local Jacobian determinant approaches zero from above:
\begin{align}
\mathcal{L}_{guard} &= \frac{1}{N} \sum_{i=1}^N \left[ \operatorname{ReLU}\bigl(J_{safe\_min} - \det(\mathbf{J}_i)\bigr) \right]^2. 
\end{align}

This strictly prevents $\det(\mathbf{J})$ from vanishing or turning negative, thus rigorously preserving the bijectivity, smoothness, and topological integrity of the mapping over the entire domain.

\subsubsection{AGM dynamic evolution callback mechanism}

To efficiently realize the closed-loop iteration of physical gradient capture, geometric grid deformation, coordinate update, and manifold reconstruction during training, a unified callback mechanism is designed. 

In Stage 1, the AGM callback executes the following closed loop every 100 training steps:
\begin{itemize}
\item \textbf{Gradient field capture and forward decoupling:} Pause back-propagation of the governing equation residual, pass the current collocation point batch through the physical backbone network, activate the first-order automatic differentiation graph, capture the latest spatial gradient magnitudes, detach them from the physical optimization graph, and convert them to static driving sources.
\item \textbf{Monitoring operator reconstruction and manifold update:} Substitute the extracted gradient field into the governing equations to compute the dynamic monitoring indicator \(\omega(X)\), calculate the Jacobian determinant of the current transformation, and trigger an independent variational optimization step for the mapping network parameters \(\Theta_{MAP}\) within the callback.
\item \textbf{Coordinate mapping dynamic overwrite:} As \(\Theta_{MAP}\) updates, the manifold coordinates \(\xi\) output by the mapping network are refreshed in real time and passed as new input features to the subsequent FFM, reconstructing a smoother manifold space more amenable to spectral approximation for the next 100 steps of physical field solution.
\end{itemize}

The mapping network does not participate in PDE loss back-propagation of the backbone network; its parameter updates are independently driven by the AGM callback mechanism.

\subsection{Gaussian Fourier feature mapping (FFM)}
\label{subsec:ffm}

NTK theory indicates that deep fully connected networks exhibit exponentially decaying convergence rates for high-frequency features when approximating complex physical fields with multi-scale, high-wavenumber components—the spectral bias phenomenon. To achieve "pseudo-low-frequency" mapping of flow field features, GAC-PINN introduces a Gaussian Fourier feature embedding layer \cite{tancik2020fourier} before the transformed manifold coordinates enter the backbone network. Let the transformed spatiotemporal manifold coordinate vector be \(z_{MAP} = [\xi, t]^T \in \mathbb{R}^d\). The high-dimensional spectral feature projection is defined as:
\begin{equation}
\gamma(z_{MAP}) = \begin{bmatrix}
\sin(2\pi B z_{MAP}) \\
\cos(2\pi B z_{MAP})
\end{bmatrix} \in \mathbb{R}^{2m},
\end{equation}
where \(m = 128\) is the half-dimension of the Fourier feature mapping, yielding an embedded feature dimension of \(2m = 256\). The projection matrix \(B \in \mathbb{R}^{m \times d}\) is generated once during initialization and fixed as a static constant, with elements sampled from an isotropic multivariate Gaussian distribution:
\begin{equation}
B_{jk} \sim \mathcal{N}(0, \delta_{scale}^2),
\end{equation}
with Gaussian standard deviation \(\delta_{scale} = 5.0\). By explicitly mapping the low-dimensional spatiotemporal domain to the high-dimensional spectral space spanned by the Gaussian prior basis, FFM pre-stretches the frequency response distribution of the input signal, equipping the physical backbone network with the capability to capture high-wavenumber information of strongly nonlinear physical interfaces from the early stages of training.

\subsection{Adaptive hard-constraint network}
\label{subsec:ahc}
The hard-constraint ansatz has been widely applied in the field of PINN research to eliminate boundary loss and prevent weight collapse. This work employs a variant with a learnable bandwidth and a stop-gradient operator, which is an engineering enhancement intended to stabilize high-order differentiation.
\subsubsection{Gradient-truncated hard constraints for non-periodic boundary conditions}

For non-periodic boundaries (e.g., homogeneous Dirichlet boundaries for 2D Poisson, fixed-value boundaries for 1D Allen-Cahn), the general hard-constraint ansatz is defined as \cite{straub2025hardconstraining,lu2021physics}:
\begin{equation}
u(X) = g(X) + t \cdot \mathcal{D}_{space}(X; \kappa(X)) \cdot \mathcal{N}_{phy}(\gamma(z_{MAP}); \Theta_{PHY}),
\end{equation}
where \(g(X)\) is the basis function exactly satisfying initial conditions, \(\mathcal{N}_{phy}\) is the physical backbone network, \(\gamma(\cdot)\) is the Fourier feature embedding mapping, and \(z_{MAP}\) is the transformed manifold coordinate from AGM. The spatial boundary distance factor \(\mathcal{D}_{space}(X)\) is constructed as:
\begin{equation}
\mathcal{D}_{space}(X) = 1.0 - \exp\left(-\lfloor \kappa_{adaptive}(X) \rfloor_{sg} \cdot \mathcal{D}_{mesh}(X)\right),
\end{equation}
where \(\mathcal{D}_{mesh}(X)\) is the distance from the physical space coordinate to the boundary. The bandwidth coefficient \(\kappa_{adaptive}\) is predicted in real time by a lightweight gating network:
\begin{equation}
\kappa_{adaptive}(X) = \kappa_{min} + (\kappa_{max} - \kappa_{min}) \cdot \operatorname{Sigmoid}\left(\mathcal{N}_{gate}(\nu_{gate}; \Phi_{gate})\right), 
\end{equation}
with input vector \(\nu_{gate} = [x, t]\) comprising spatial and temporal coordinates; \(\Phi_{gate}\) denotes the trainable gating network parameters, and \(\lfloor\cdot\rfloor_{sg}\) denotes the stop-gradient operator.

When computing high-order spatial derivatives in the PDE operator \(\mathcal{P}[u]\) (e.g., \(\partial^2 u/\partial x^2\) or \(\partial^4 u/\partial x^4\)) via automatic differentiation, \(\lfloor\kappa_{adaptive}(X)\rfloor_{sg}\) is treated as a spatial constant independent of network parameters. This design cuts the gradient propagation path from \(\Phi_{gate}\) to the high-order residual automatic differentiation graph, mathematically eliminating the potential influence of the bandwidth prediction network on the numerical stiffness of high-order PDE derivative terms.

\subsubsection{Intrinsic periodic construction for periodic boundary conditions}

For the periodic boundary conditions $u(-1, t) = u(1, t)$ of the 1D Burgers equation, our framework employs a fully periodic harmonic basis substitution rather than algebraic truncation terms. The periodic ansatz is defined as:
\begin{equation}
u(X) = e^{-\lambda t} \cdot g(x) + \left(1.0 - e^{-\lambda t}\right) \cdot \mathcal{N}_{phy}\left(\sin(\pi x), \cos(\pi x), \gamma_t(t); \Theta_{PHY}\right),
\end{equation}
where $g(x) = -\sin(\pi x)$ exactly matches the initial condition $u(x, 0) = -\sin(\pi x)$. Since the network input layer is explicitly constructed as $\sin(\pi x)$ and $\cos(\pi x)$ basis forms, the output intrinsically satisfies full periodic continuity $u(x+2, t) \equiv u(x, t)$ regardless of weight evolution. The temporal decay factor $e^{-\lambda t}$ ensures exact recovery of initial conditions at $t=0$, while $1.0 - e^{-\lambda t}$ progressively delegates solution evolution to the backbone network as $t$ increases, simultaneously imposing initial and periodic boundary conditions without additional penalty terms.

\subsection{Three-stage training strategy}

To ensure stable convergence of the relative \(L^2\) error to high precision for strongly nonlinear evolution equations, we adopt a three-stage training strategy:

\begin{itemize}
\item \textbf{Stage 1: Adam with AGM cooperative evolution:} Physical network weights \(\Theta_{PHY}\) and mapping network weights \(\Theta_{MAP}\) evolve asynchronously in a dual-track manner. AGM performs offline forward updates via the callback mechanism every 100 training steps. During this stage, the system conducts global topological exploration over the large-scale spatiotemporal domain; high-gradient regions are initially localized and stretched on the geometric manifold, establishing the fundamental solution topology.
\item \textbf{Residual-based adaptive refinement (RAR) \cite{mao2023physicsinformed}:} After Stage 1, one RAR step is performed: 60,000 candidate points are sampled from the full domain, PDE residual magnitudes are computed, and the 2,500 points with the largest residuals are appended to the training set, focusing computational resources on regions with maximal physical residuals.
\item \textbf{Stage 2: Freeze manifold, fine-tune physical layers:} All trainable parameters \(\Theta_{MAP}\) of AGM are frozen, eliminating potential high-frequency oscillations from coordinate transformation grids in late training and fixing the manifold geometry. Computational resources are concentrated on resolving high-gradient regions captured by RAR using Adam fine-tuning of \(\Theta_{PHY}\).
\item \textbf{Stage 3: L-BFGS final convergence, $N_3 = 3{,}000$ iterations}: The L-BFGS second-order optimizer is employed for final convergence with full-batch computation, leveraging curvature information to achieve precise local minimization of $\mathcal{L}_{PDE}$. Only $\{\Theta_{PHY}, \Phi_{gate}\}$ are updated; $\Theta_{MAP}$ remains frozen to preserve the established manifold geometry.

\end{itemize}

\subsection{Algorithm}

The complete GAC-PINN training procedure is summarized in Algorithm 1.

\begin{algorithm}[H]
\SetAlgoLined
\caption{GAC-PINN Training Algorithm}  
\label{alg:gacpinn}

\KwIn{PDE operator \(\mathcal{P}\), computational domain \(\Omega\) and time interval \([0,T]\), boundary and initial conditions, total collocation points \(\mathcal{N}_{pde}\), hyperparameters \(\lambda_{equi}, \lambda_{barrier}, J_{safe\_min}\)}
\KwOut{Solution space \(u \in U\), neural network parameters \(\Theta_{NN}\)}
Initialize physical network \(\mathcal{N}_{phy}\) (parameters \(\Theta_{phy}\)), AGM network \(\mathcal{N}_{AGM}\) (parameters \(\Theta_{AGM}\)), bandwidth gating network \(\mathcal{N}_{gate}\) (parameters \(\Phi_{gate}\));\\
Generate Fourier projection matrix \(B\) by fixed sampling \(B_{jk} \sim \mathcal{N}(0, \delta_{scale}^2)\), construct FFM embedding \(\gamma(\cdot)\);\\
Select hard-constraint ansatz according to boundary topology, construct output transformation \(\mathcal{T}\), bind into complete forward propagation chain:\\
\quad \(u_{NN}(X) = \mathcal{T}(\mathcal{N}_{phy}(\gamma(\mathcal{N}_{AGM}(X; \Theta_{AGM})); \Theta_{phy}); \Phi_{gate})\);\\
\textbf{Stage 1} (cooperative training, \(N_1 = 12000\)):\\
\For{\(iter = 1\) to \(N_1\)}{
    fix \(\Theta_{AGM}\), update \(\Theta_{phy}, \Phi_{gate}\) with Adam minimizing \(\mathcal{L}_{PDE} = \frac{1}{\mathcal{N}_{pde}}\sum_{i=1}^{\mathcal{N}_{pde}} \| \mathcal{P}[u_{NN}](X_i) \|_2^2\);\\
    \If{\(iter \bmod 100 = 0\)}{
        fix \(\Theta_{phy}, \Phi_{gate}\), update \(\Theta_{AGM}\) with Adam minimizing \(\mathcal{L}_{AGM} = 2.0 \cdot \mathcal{L}_{equi} + \lambda_{barrier} \cdot \mathcal{L}_{guard}\);
    }
}
Execute RAR: select \(N_{refine} = 2500\) points with largest residuals from \(N_{cand} = 60000\) candidate points, append to training set;\\
\textbf{Stage 2} (frozen mapping fine-tuning, \(N_2 = 5000\)): freeze \(\Theta_{AGM}\), update \(\Theta_{phy}, \Phi_{gate}\) with Adam (lr \(= 10^{-4}\)) minimizing \(\mathcal{L}_{PDE}\);\\
\textbf{Stage 3} (L-BFGS refinement, \(N_3 = 3000\)): update all parameters \(\{\Theta_{phy}, \Phi_{gate}\}\) with L-BFGS minimizing \(\mathcal{L}_{PDE}\);\\
\Return trained GAC-PINN model.
\end{algorithm}

\section{Experimental setup}

\subsection{Computational environment}
All experiments were developed and tested on a Linux computing platform equipped with an NVIDIA Tesla V100 GPU, utilizing PyTorch and the DeepXDE scientific computing library. The benchmark datasets and problem configurations for the Burgers equation, Allen-Cahn equation, and two-dimensional Poisson equation were adopted from the open-source code repository associated with the gradient-enhanced physics-informed neural networks (gPINN) work by Yu et al.~\cite{yu2022gradienenhanced}, which is publicly available at \url{https://github.com/lu-group/gpinn}. For the two‑dimensional unsteady cylinder‑wake benchmark, we adopt the governing‑equation configuration from PINNsFormer~\cite{zhao2024pinnsformer}, while the reference ground‑truth flow‑field data originate from the Nek5000 simulation database of Raissi et al.~\cite{raissi2019physicsinformed}. To eliminate the influence of random perturbations on optimization trajectories, all experiments were performed with five independent random seeds, and results are reported as mean $\pm$ standard deviation. For clarity, the error distribution figures and the specific error values reported in Section~\ref{sec:num_exp} are obtained from a single representative run with the same fixed random seed.

\subsection{Benchmark problems and baselines}

\subsubsection{PDE benchmark problems}

\begin{itemize}
\item \textbf{1D Burgers equation:} The kinematic viscosity is set to \(\nu = 0.01/\pi\), testing the model's capability for localized capture of nonlinear fluid shock fronts under convection-dominated conditions.
\item \textbf{1D Allen-Cahn equation:} This equation includes a high-order nonlinear cubic reaction source term (\(\epsilon = 0.001\)) with sharp phase interface rotation and spatiotemporal evolution, evaluating the efficacy of the time-adaptive hard constraint.
\item \textbf{2D Poisson equation:} The method of manufactured solutions introduces a high-order exponential parameter (\(a = 10\)), constructing an extremely sharp spatial singularity peak at the domain center to assess the manifold adaptive compression performance of the model for complex static spatial gradients.
\item \textbf{2D Navier-Stokes equation:} The unsteady flow past a cylinder is considered with parameters $\lambda_1 = 1$ and $\lambda_2 = 0.01$, serving as a benchmark to evaluate the model's generalization capability for complex fluid dynamics and pressure field reconstruction.
\end{itemize}

\subsubsection{Baseline models}

Five representative baseline models spanning the evolutionary spectrum from traditional penalty methods to recent dynamic weighting and adaptive sampling strategies were selected for comparison:

\begin{itemize}
\item \textbf{Vanilla PINN} (Raissi et al., 2019): Standard fully connected architecture with mean squared error soft constraints for boundary conditions.
\item \textbf{GPINN} (Yu et al., 2022): Explicitly embeds first-order spatial derivatives of PDE residuals as gradient supervision in the total loss.
\item \textbf{Causal PINN} (Wang et al., 2024): Explicit temporal weighting scheme based on forward-moving temporal residuals.
\item \textbf{RAR-PINN} (Mao \& Meng, 2023): Discrete incremental greedy resampling strategy dynamically appending high-residual collocation points.
\item \textbf{RAMS-PINN} (Ouyang et al., 2026): Dynamically modulates loss term weights or multi-scale basis functions using real-time residual magnitude evolution.
\end{itemize}

Moreover, for the 2D Navier-Stokes equation, which involves strong convective unsteadiness and vortex shedding dynamics that differ substantially from the preceding benchmark problems, we adopt a separate set of baselines specifically designed for fluid mechanics applications. These include PINNsformer \cite{zhao2024pinnsformer}, a transformer-based architecture for PDE solving; DD-PINN \cite{khademi2024ddpinn}, which employs a domain-decomposition strategy; and TSA-PINN \cite{khademi2025tsapinn}, which incorporates trainable sinusoidal activation functions. This targeted selection ensures physically meaningful comparisons tailored to Navier-Stokes problem. 

Due to the unavailability of open-source code for DD-PINN and the irreproducibility of PINNsformer's reported results in our environment, we directly cite their published numerical results under the same initial and boundary conditions, sampling strategies, and relative $L_2$ error metrics, rendering cross-paper comparisons objectively valid.

\subsection{Evaluation metrics}

\subsubsection{Relative \(L^2\) error}

The relative \(L^2\) error is adopted as the primary accuracy metric:
\begin{equation}
\epsilon_{L^2} = \frac{\| u_{pred} - u_{true} \|_2}{\| u_{true} \|_2} = \sqrt{\frac{\sum_{i=1}^N |u_{pred}(X_i) - u_{true}(X_i)|^2}{\sum_{i=1}^N |u_{true}(X_i)|^2}}.
\end{equation}

\subsubsection{Convergence order}

To quantify the error decay with mesh refinement, we define the convergence order \(p\) via the power-law relationship:
\begin{equation}
E = C \cdot N_{pde}^{-p},
\end{equation}
where \(E\) is the relative \(L^2\) error and \(C\) is a constant. Taking logarithms gives \(\log E = \log C - p\log N_{pde}\); thus \(p\) is estimated as the negative slope of the linear fit in the log-log plane. For each seed, an individual \(p_i\) is obtained by fitting its errors across all \(N_{pde}\). A global \(p_{\text{global}}\) is similarly estimated from the geometric mean errors at each \(N_{pde}\).

\section{Numerical experiments}
\label{sec:num_exp}
\subsection{Burgers equation: convective shock flow field}

The Burgers equation serves as a simplified nonlinear form of the Navier-Stokes equations in fluid mechanics, commonly used to model convection-dominated flow and shock formation. The governing equation is:
\begin{equation}
\frac{\partial u}{\partial t} + u\frac{\partial u}{\partial x} = \nu \frac{\partial^2 u}{\partial x^2}, \quad x \in [-1, 1], \ t \in [0, 1],
\end{equation}
with \(\nu = 0.01/\pi\), initial condition \(u(x, 0) = -\sin(\pi x)\), and periodic boundary conditions \(u(-1, t) = u(1, t)\).

As \(t \rightarrow 1.0\), the nonlinear convection term drives rapid steepening of the solution field near \(x = 0\), evolving an extremely thin fluid shock front. 

To demonstrate the dynamic adaptability of the AGM framework, Fig.~\ref{fig:agm_burgers} illustrates the mesh redistribution on the 1D Burgers equation across snapshots $t = 0.2, 0.5,$ and $0.8$. The left panels show the coordinate mapping ($\xi$ to $x$) with collocation nodes, while the right panels display the corresponding field solutions $u(x,t)$ and high-gradient shock regions. Fig.~\ref{fig:agm_burgers} reveals the core regulatory features of AGM:

\begingroup
\centering
\includegraphics[width=0.7\textwidth]{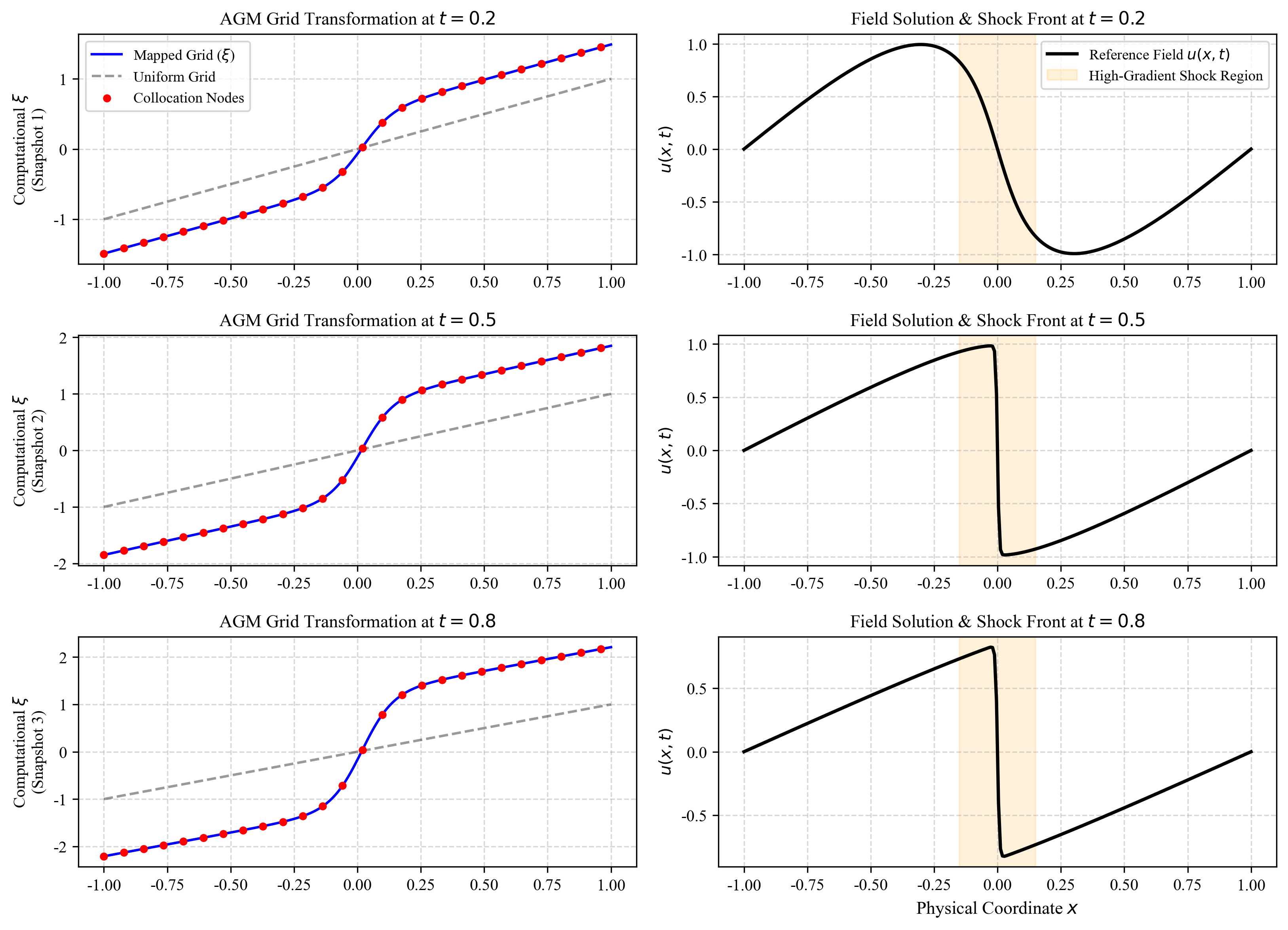}
\captionof{figure}{Adaptive grid evolution and physical field solutions for the 1D Burgers' equation across representative time snapshots ($t = 0.2, 0.5, 0.8$).}
\label{fig:agm_burgers}
\endgroup
\medskip

\begin{itemize}
\item \textbf{Shock-tracking and dynamic concentration:} 
As shown in the left-hand panels, as time evolves from $t = 0.2$ to $t = 0.8$ and the shock front of the physical field on the right becomes increasingly steep, the mapped grid ($\xi$) and collocation nodes (red dots) exhibit a significant and monotonically intensifying S-shaped deformation in the central high-gradient region ($x \in [-0.15, 0.15]$). Meanwhile, the range of the vertical computational domain coordinate dynamically expands from approximately $\pm 1.5$ to around $\pm 2$, intuitively reflecting the AGM framework's continuously enhanced local compression and node-refinement capabilities in the high-gradient region over time.

\item \textbf{Precise field-grid correspondence:} 
Comparing the physical field solution $u(x,t)$ on the right with the coordinate mapping on the left indicates that the AGM module can perceive drastic local gradient changes in real-time, precisely allocating dense computational nodes in the core shock region where physical variations are most severe, thereby achieving efficient adaptive deployment of computational resources.

\item \textbf{Boundary and smooth domain preservation:} 
In the low-gradient regions near the physical boundaries ($x = \pm 1$), the grid transformation smoothly transitions and approaches a linear uniform distribution without causing unnecessary distortion, thereby ensuring the numerical stability of boundary condition enforcement.

\item \textbf{Optimization compatibility for stiff problems:} 
While maintaining topological homeomorphy, continuous differentiability, and a constant total number of collocation points, this method avoids the loss-function discontinuities brought by traditional dynamic remeshing techniques (such as RAR), ensures the smoothness of the variational optimization landscape, and thus perfectly matches the high-precision convergence requirements of second-order optimizers (such as L-BFGS).
\end{itemize}

GAC-PINN demonstrates superior adaptive capture capability for the Burgers shock problem. The FFM, through nonlinear stretching via Gaussian random projection, enables the backbone network to establish strong representations of high-wavenumber shock fronts from early optimization stages. AGM captures strong spatial physical gradients at the shock front in real time, driving spontaneous local high-density aggregation of computational manifold nodes at \(x = 0\) through the global synergy of the equidistribution loss \(\mathcal{L}_{equi}\). RAR appends 2,500 high-residual points in the shock layer as a supplementary refinement step; its individual contribution is examined in the ablation study. The GAC-PINN predicted flow field aligns closely with the analytical solution, achieving a global relative \(L^2\) error of \(6.785\times10^{-5}\) and effectively eliminating non-physical numerical oscillations.

Fig.~\ref{fig:burgers_error} presents the absolute error distribution between the predicted and reference solutions across the full spatio-temporal domain.

\begingroup
\centering
\includegraphics[width=0.8\textwidth]{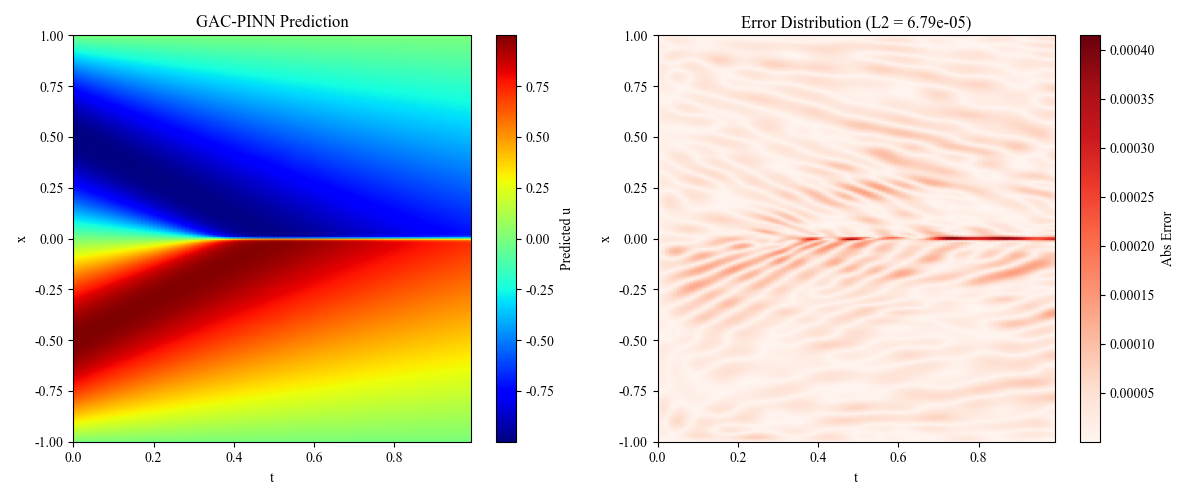}
\captionof{figure}{Absolute error distribution of GAC-PINN for the Burgers equation.}
\label{fig:burgers_error}
\endgroup
\medskip

From the error contours, three observations can be made: (1) the maximum error is strictly localized near the shock front trajectory around \(x \approx 0\) and \(t \rightarrow 1.0\), corresponding to the infinite gradient discontinuity in the Burgers equation; (2) no non-physical oscillations are observed on either side of the shock, benefiting from FFM's NTK spectrum reshaping and AGM's geometric widening effect; (3) the error contours vary continuously in the spatio-temporal domain without abrupt changes or discontinuities, confirming that the three-stage training strategy effectively avoids the loss landscape discontinuity caused by discrete resampling.

\subsection{2D Poisson equation: extreme gradient potential reconstruction}

The 2D Poisson equation is a canonical elliptic PDE model widely applicable to electrostatics, heat conduction, potential flow, and elasticity. The governing equation is:
\begin{equation}
\nabla^2 u(x, y) = f(x, y), \quad (x, y) \in [0, 1]^2,
\end{equation}
with homogeneous Dirichlet boundary conditions \(u|_{\partial\Omega} = 0\). The method of manufactured solutions specifies the reference solution:
\begin{equation}
u_{true}(x, y) = [16xy(1-x)(1-y)]^a, \quad a = 10,
\end{equation}
with the source term \(f(x, y)\) analytically derived from \(f = -\nabla^2 u_{true}\).

In the Poisson experiment, the operator-aware router determines \(\chi_{\mathcal{P}} = 0\) since the control operator \(\mathcal{P} = \nabla^2 u - f\) contains no temporal derivatives. The framework therefore selects the pure spatial distance function as the hard-constraint construction, dedicating all computational resources to the spatial adaptive reconstruction module. AGM, through minimization of the manifold volume dispersion \(\mathcal{L}_{equi}\), drives nonlinear shear and compression of physical space coordinates. At the central high-gradient peak region, the Jacobian determinant $\det(\mathbf{J})$ increases substantially, indicating strong local densification of collocation points. The safety barrier $J_{safe\_min} = 0.02$ prevents $\det(\mathbf{J})$ from approaching zero, thereby guaranteeing that the mapping remains bijective throughout training. Coupled with FFM and adaptive hard constraints, GAC-PINN perfectly reconstructs the sharp 2D potential field with a global relative \(L^2\) error of only \(3.077\times10^{-5}\).

To intuitively demonstrate the spatial modeling fidelity of the proposed framework, Fig.~\ref{fig:possion_3d} presents the three-dimensional surface plot of the predicted potential field $u(x, y)$. As illustrated, the network accurately captures the sharp, centralized bell-shaped profile induced by the high-order manufactured solution ($a=10$), while exhibiting exceptional smoothness and complete compliance with the homogeneous Dirichlet boundary conditions across all perimeter edges of the computational domain $[0, 1]^2$.

\begingroup
\centering
\includegraphics[width=0.55\textwidth]{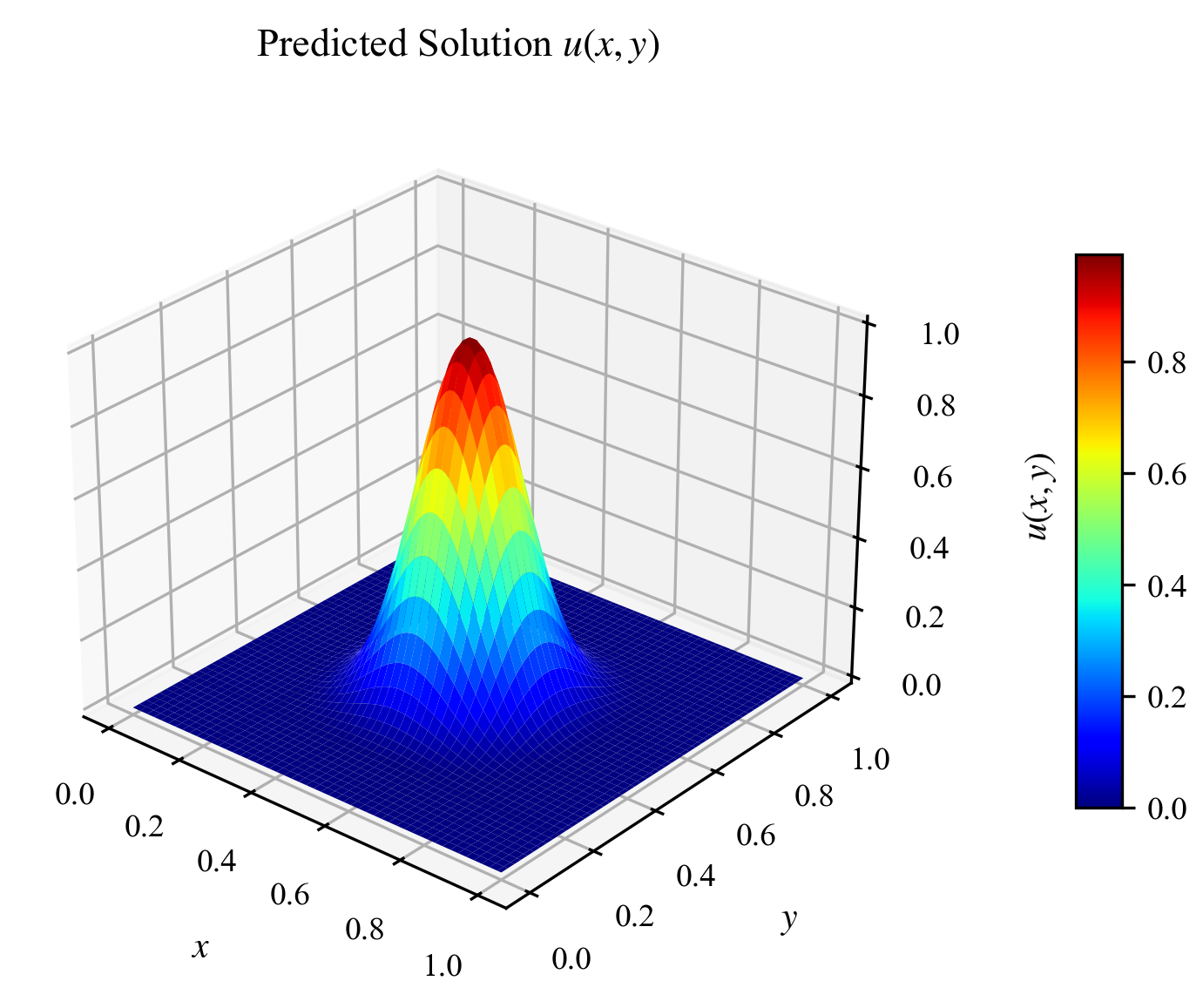}
\captionof{figure}{Three-dimensional surface plot of the GAC-PINN predicted solution $u(x, y)$ for the 2D Poisson equation.}
\label{fig:possion_3d}
\endgroup
\medskip

\begingroup
\centering
\includegraphics[width=0.9\textwidth]{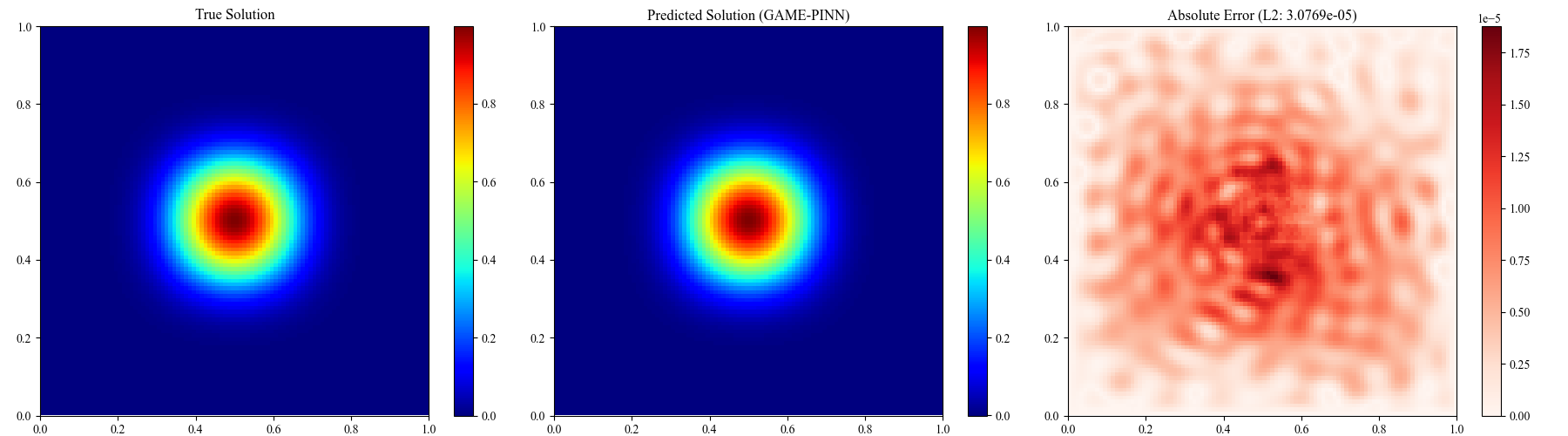}
\captionof{figure}{Comprehensive comparison for the 2D Poisson equation: (left) analytical true solution, (middle) GAC-PINN prediction, and (right) absolute error distribution field, achieving a global relative $L^2$ error of $3.077 \times 10^{-5}$.}
\label{fig:possion_comparison}
\endgroup
\medskip

To further quantify and validate the reconstruction performance, Fig.~\ref{fig:possion_comparison} provides a comprehensive three-way visual comparison comprising the analytical true solution, the GAC-PINN predicted solution, and the spatial absolute error distribution. The qualitative comparison between the true and predicted solutions confirms a near-perfect visual alignment. Furthermore, the absolute error distribution field reveals that the minor residuals are strictly localized around the steep central peak region characterized by high local curvature and gradients, whereas the extensive outer low-gradient domains maintain near-zero error levels. This outcome confirms the effectiveness of the AGM-driven adaptive collocation point aggregation strategy in tackling extreme gradient challenges.

\subsection{Allen-Cahn equation: sharp phase‑interface flow field}

The Allen-Cahn equation is a canonical nonlinear parabolic PDE describing interface evolution and reaction-diffusion processes in multi-phase materials science:
\begin{equation}
\frac{\partial u}{\partial t} - \epsilon \frac{\partial^2 u}{\partial x^2} + 5(u^3 - u) = 0, \quad x \in [-1, 1], \ t \in [0, 1],
\end{equation}
with phase interface diffusion bandwidth coefficient \(\epsilon = 0.001\), initial condition \(u(x, 0) = x^2 \cos(\pi x)\), and boundary conditions \(u(-1, t) = u(1, t) = -1\).

For the Allen-Cahn problem, the operator-aware router determines \(\chi_{\mathcal{P}} = 1\) because the governing equation contains a temporal derivative. The framework therefore selects the time-adaptive bandwidth hard constraint, which gradually releases the initial condition as \(t\) increases. Meanwhile, AGM applies its nonlinear transformation only to the spatial axis, preserving temporal independence. The global relative \(L^2\) error remains stable at \(8.996\times10^{-3}\), outperforming existing baseline models. Fig.~\ref{fig:allen_error} presents the error distribution of GAC-PINN for the Allen-Cahn sharp phase-field evolution.

\begingroup
\centering
\includegraphics[width=0.8\textwidth]{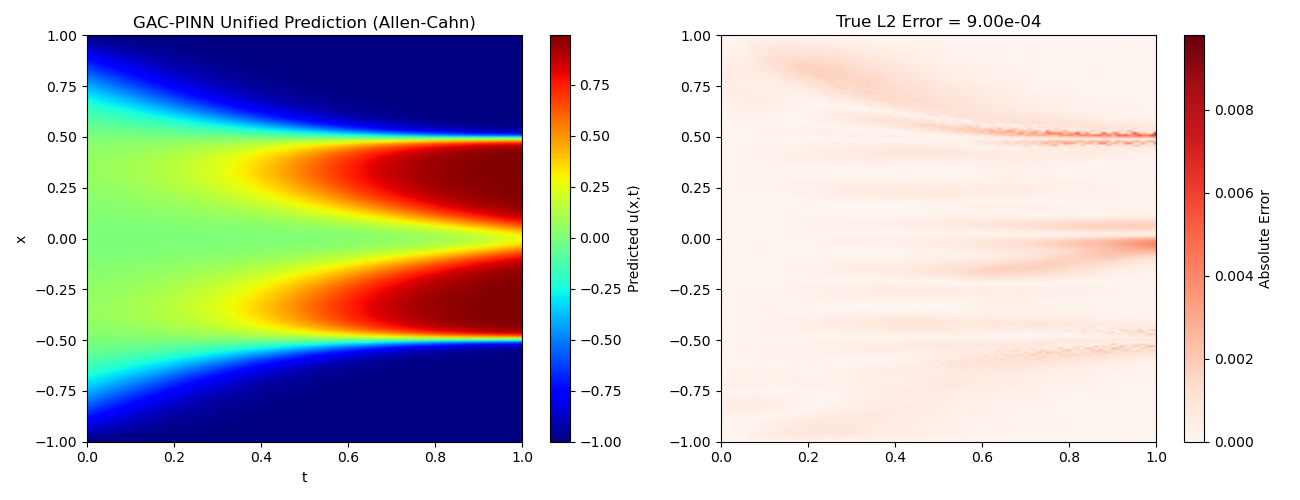}
\captionof{figure}{Absolute error distribution of GAC-PINN for the Allen-Cahn equation.}
\label{fig:allen_error}
\endgroup
\medskip

The left panel of Fig.~\ref{fig:allen_error} depicts the dynamic evolution from the initial state to multi-phase separation. The right panel shows that the vast majority of the spatio-temporal domain exhibits extremely low absolute error (approaching \(10^{-4}\)); minor local error elevations are concentrated at steep phase interface gradient regions and near \(t = 1.0\), consistent with the numerical challenges typically encountered in deep learning-based PDE solutions.

\subsection{2D Navier‑Stokes equation: incompressible flow field}
The 2D Navier‑Stokes equation is a set of parabolic partial differential equations describing incompressible fluid dynamics. As fundamental governing equations in fluid mechanics, they are widely adopted in scientific investigations and engineering applications to model the motions of fluid media such as water and air. The governing equations are as follows:
\begin{equation}
\begin{aligned}
& \frac{\partial u}{\partial t}+\lambda_{1}\left(u \frac{\partial u}{\partial x}+v \frac{\partial u}{\partial y}\right)=-\frac{\partial p}{\partial x}+\lambda_{2}\left(\frac{\partial^{2} u}{\partial x^{2}}+\frac{\partial^{2} u}{\partial y^{2}}\right), \\
& \frac{\partial v}{\partial t}+\lambda_{1}\left(u \frac{\partial v}{\partial x}+v \frac{\partial v}{\partial y}\right)=-\frac{\partial p}{\partial y}+\lambda_{2}\left(\frac{\partial^{2} v}{\partial x^{2}}+\frac{\partial^{2} v}{\partial y^{2}}\right),
\end{aligned}
\label{eq:2d_ns}
\end{equation}
where $u(t,x,y)$ and $v(t,x,y)$ stand for velocity components in the $x$ and $y$ directions, and $p(t,x,y)$ is fluid pressure. In this work, the coefficients are set to $\lambda_1=1$ and $\lambda_2=0.01$. Table~\ref{tab:ns_error} summarizes the relative $L_2$ errors of all compared models for the three output fields (u, v, and p).

As shown in Table~\ref{tab:ns_error}, Vanilla PINN achieves errors of $7.95\times10^{-3}$ for u and for $2.42\times10^{-2}$ v in the velocity field, but its pressure field error reaches as high as 17.27, indicating that conventional PINNs struggle to accommodate the wide amplitude span of pressure gradients in the absence of explicit boundary constraints, leading to severe degradation in pressure reconstruction. DD-PINN reduces the pressure error to $0.28$ through a domain-decomposition strategy, while also improving the velocity errors over the vanilla PINN; PINNsformer similarly achieves a pressure error of $0.28$. The velocity errors of TSA-PINN are on the order of $10^{-2}$, while those of our GAC-PINN attain a comparable order of magnitude with slight improvements, thereby confirming its effective reconstruction of the velocity field. Meanwhile, GAC-PINN drastically reduces the pressure error to $2.66\times10^{-2}$, representing a nearly one-order-of-magnitude reduction from the $0.28$ achieved by DD-PINN and PINNsformer. This result fully demonstrates its capability for simultaneous accurate reconstruction of both velocity and pressure fields in unsteady flows, as well as its favorable generalization performance.

\begin{center}
\centering
\captionof{table}{Relative $L_2$ errors of different models on 2D Navier‑Stokes equation.}
\label{tab:ns_error}
\begin{tabular}{cccc}
\toprule
Models & Relative L2 Error ($u$) & Relative L2 Error ($v$) & Relative L2 Error ($p$) \\
\midrule
Vanilla PINN & $7.95\times 10^{-3}$ & $2.42\times 10^{-2}$ & $17.27$ \\
PINNsformer & -- & -- & $0.28$ \\
DD‑PINN & $6.2\times 10^{-3}$ & $1.9\times 10^{-2}$ & $0.28$ \\
TSA‑PINN & $2.94\times 10^{-2}$ & $5.44\times 10^{-2}$ & -- \\
GAC‑PINN & $5.47\times 10^{-3}$ & $1.62\times 10^{-2}$ & $2.66\times 10^{-2}$ \\
\bottomrule
\end{tabular}
\end{center}
\medskip

This improvement is primarily attributed to two factors: the gradient-adaptive balancing strategy, which dynamically adjusts the weights of individual loss terms and effectively prevents the pressure gradient from being dominated by velocity residuals during back-propagation; and the physics-constrained attention mechanism, which enhances the network's capacity to capture localized sharp features in the pressure field.

Figure~\ref{fig:ns_error} plots the exact pressure field \(p(x,y)\), the GAC‑PINN prediction, and the absolute error at \(t=10.0\). It can be observed that the predicted pressure contour faithfully reproduces the main spatial patterns, extreme-value regions and sharp gradient structures of the reference solution. The absolute error map reveals that large errors mainly concentrate near steep pressure gradients, whereas the majority of the computational domain maintains a low error level.

\begingroup
\centering
\includegraphics[width=1.0\textwidth]{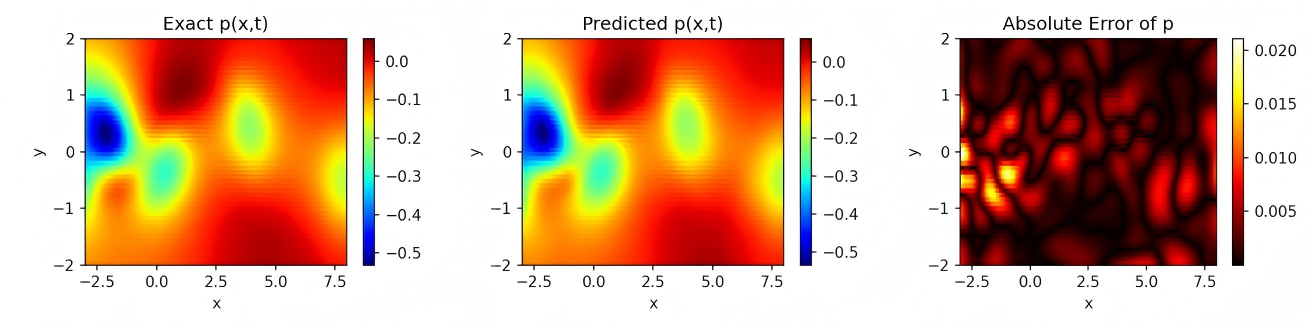}
\captionof{figure}{Pressure field reconstruction results for the 2D Navier‑Stokes equation.}
\label{fig:ns_error}
\endgroup
\medskip

The contour comparison further verifies the quantitative error listed in Table~\ref{tab:ns_error}. Given the relative $L_2$ error of $2.66\times10^{-2}$ for pressure, the contour plots and error distributions jointly demonstrate that GAC‑PINN is capable of capturing complex localized flow structures and achieving high‑fidelity pressure field reconstruction for incompressible Navier‑Stokes problems.

\section{Discussion}

\subsection{Prediction accuracy comparison}

Table 1 presents the relative \(L^2\) errors of all models across the three benchmark problems. GAC-PINN achieves optimal accuracy in all tests, with particularly pronounced advantages on steep gradient or discontinuous problems.

\begin{center}
\captionof{table}{Relative \(L^2\) error comparison (mean \(\pm\) standard deviation over five independent runs).}
\label{tab:accuracy}
\begin{tabular}{lccc}
\toprule
Models & Burgers & Poisson & Allen-Cahn \\
\midrule
Vanilla PINN  & \((2.696\pm1.426)\times10^{-2}\) & \((5.269\pm0.877)\times10^{-4}\) & \((5.145\pm4.316)\times10^{-2}\) \\
RAR-PINN      & \((6.389\pm9.970)\times10^{-3}\) & \((5.313\pm0.864)\times10^{-4}\) & \((3.023\pm0.902)\times10^{-2}\) \\
GPINN         & \((8.512\pm1.680)\times10^{-3}\) & \((1.065\pm0.588)\times10^{-3}\) & \((1.951\pm0.308)\times10^{-2}\) \\
Causal PINN   & \((1.305\pm0.315)\times10^{-2}\) & \((2.914\pm0.915)\times10^{-3}\) & \((2.527\pm1.064)\times10^{-3}\) \\
RAMS-PINN     & \((1.082\pm1.231)\times10^{-2}\) & \((2.242\pm0.420)\times10^{-3}\) & \((7.689\pm0.013)\times10^{-1}\) \\
\bfseries GAC-PINN & \(\mathbf{(1.747\pm0.450)\times10^{-4}}\) & \(\mathbf{(2.868\pm0.947)\times10^{-5}}\) & \(\mathbf{(1.756\pm0.712)\times10^{-3}}\) \\
\bottomrule
\multicolumn{4}{l}{\footnotesize $^{\dagger}$ Values are reported as mean \(\pm\) std over five random seeds (1234, 5678, 9012, 3456, 7890).} \\
\end{tabular}
\end{center}

For the low-viscosity Burgers equation, GAC-PINN achieves a relative \(L^2\) error of \((1.747\pm0.450)\times10^{-4}\), more than one order of magnitude lower than the optimal baseline RAMS-PINN (\((1.082\pm1.231)\times10^{-2}\)) and with substantially smaller variance. For the 2D Poisson problem, GAC-PINN achieves \((2.868\pm0.947)\times10^{-5}\), more than an order of magnitude lower than the best baseline RAR-PINN (\((5.313\pm0.864)\times10^{-4}\)). For the Allen-Cahn equation, GAC-PINN attains \((1.756\pm0.712)\times10^{-3}\), comparable to Causal PINN (\((2.527\pm1.064)\times10^{-3}\)) but with noticeably smaller standard deviation, indicating more stable convergence without manual tuning of temporal weights.

\subsection{Computational efficiency analysis}

Table 2 summarizes computational cost metrics across models for the 1D Burgers problem, and Fig.~\ref{fig:efficiency} presents the evolution of relative \(L^2\) error with training iterations.

\begin{center}
\captionof{table}{Computational cost comparison.}
\label{tab:efficiency}
\begin{tabular}{lccccc}
\toprule
Model & Parameters & Peak Memory (GB) & Total Time (s) & Epochs to Converge & \(L^2\) Error \\
\midrule
Vanilla PINN  & 921    & 0.163  & 1103 & 15000  & \((2.696\pm1.426)\times10^{-2}\) \\
RAR-PINN      & 2,241  & 0.228  & 4213 & 67203  & \((6.389\pm9.970)\times10^{-3}\) \\
GPINN         & 21,057 & 0.684  & 1457 & 38308  & \((8.512\pm1.680)\times10^{-3}\) \\
Causal PINN   & 52,609 & 0.624  & 728  & 20000  & \((1.305\pm0.315)\times10^{-2}\) \\
RAMS-PINN     & 20,601 & 6.374  & 989  & 21000  & \((1.082\pm1.231)\times10^{-2}\) \\
\bfseries GAC-PINN & 104,708 & 3.776 & 1117 & 18514 & \(\mathbf{(1.747\pm0.450)\times10^{-4}}\) \\
\bottomrule
\end{tabular}
\end{center}
\medskip

\begingroup
\centering
\includegraphics[width=0.75\textwidth]{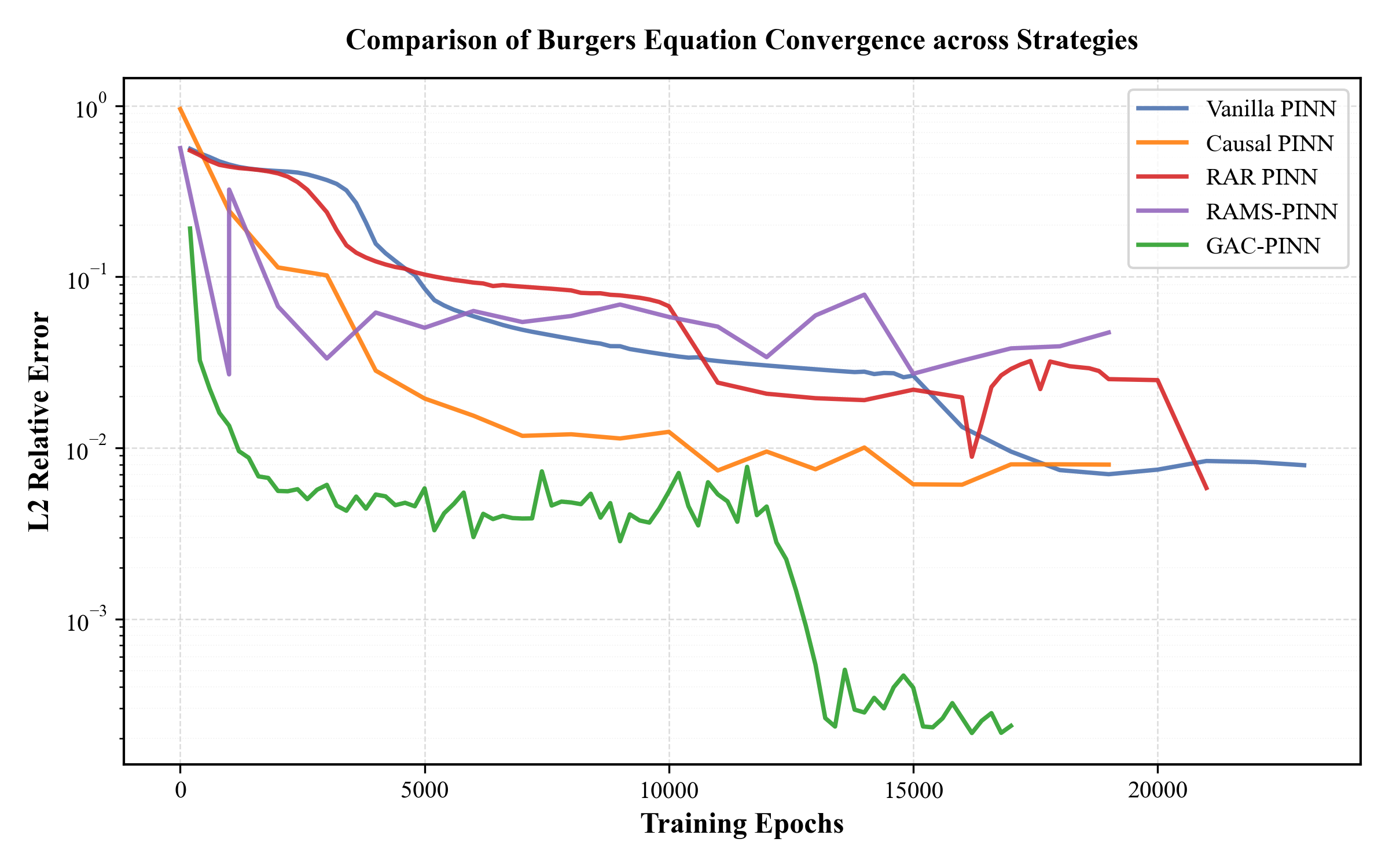}
\captionof{figure}{Evolution of relative \(L^2\) error with training iterations across all models.}
\label{fig:efficiency}
\endgroup
\medskip

GAC-PINN converges in 18,514 epochs to a relative \(L^2\) error of \((1.747\pm0.450)\times10^{-4}\), reaching the target \(10^{-4}\) precision level within approximately 12,000 steps. This convergence speed is faster than all re-implemented baselines. In terms of wall-clock time, GAC-PINN requires 1,117 s, which is comparable to Vanilla PINN (1,103 s) and substantially lower than RAR-PINN and GPINN. Although GAC-PINN has the largest parameter count and peak memory footprint, the cost increase is accompanied by a reduction of more than one order of magnitude in the relative \(L^2\) error compared with the best baseline RAMS-PINN (\((1.082\pm1.231)\times10^{-2}\)), yielding a favorable accuracy-efficiency trade-off.

\subsection{Ablation studies}

To quantitatively analyze the contribution of each core module to global convergence accuracy, seven ablation experiments were conducted on the 1D Burgers equation with a fixed random seed:

\begin{enumerate}
\item[A.] \textbf{Fixed-HC baseline:} A standard fully connected network with a fixed hard constraint.
\item[B.] \textbf{Adaptive-HC baseline:} Replaces the fixed hard constraint with an adaptive bandwidth hard constraint.
\item[C.] \textbf{B + RAR:} Adds residual-based adaptive refinement on top of B.
\item[D.] \textbf{B + AGM:} Adds adaptive grid mapping on top of B.
\item[E.] \textbf{D + RAR:} Adds RAR on top of D, combining AGM and RAR.
\item[F.] \textbf{D + FFM:} Adds Fourier feature mapping on top of D, combining AGM and FFM.
\item[G.] \textbf{GAC-PINN:} Full model with AGM, FFM, RAR, and adaptive hard constraints.
\end{enumerate}

Table~\ref{tab:ablation} summarizes the results.
\begin{center}
\small
\captionof{table}{Ablation study results for the 1D Burgers equation.}
\label{tab:ablation}
\begin{tabular}{lcccc}
\toprule
Experiment & Components & Relative \(L^2\) Error & Training Time & Parameters \\
\midrule
A: Fixed-HC base PINN      & MLP + fixed hard constraint      & \(6.890\times10^{-3}\)  & 532.07 s    & 67,843  \\
B: Adaptive-HC base PINN   & A + adaptive hard constraint     & \(1.363\times10^{-3}\)  & 487.26 s    & 67,843  \\
C: B + RAR                & B + residual-based adaptive refinement & \(3.905\times10^{-3}\) & 603.58 s & 67,843 \\
D: B + AGM                & B + adaptive grid mapping        & \(\mathbf{3.410\times10^{-4}}\) & 819.55 s & 72,196 \\
E: D + RAR                & B + AGM + RAR                    & \(1.725\times10^{-3}\)  & 743.32 s    & 72,196  \\
F: D + FFM                & B + AGM + Fourier feature mapping & \(1.572\times10^{-4}\) & 626.61 s    & 104,708 \\
G: GAC-PINN (Full)        & F + RAR                          & \(\mathbf{6.785\times10^{-5}}\) & 1,173.30 s & 104,708 \\
\bottomrule
\multicolumn{5}{l}{\footnotesize All results are obtained with a fixed random seed (4321) under identical settings.} \\
\end{tabular}
\end{center}

The ablation results reveal several key insights. Firstly, the adaptive bandwidth hard constraint reduces the relative \(L^2\) error considerably compared with the fixed-HC baseline, confirming the value of spatially-varying boundary transition widths. Secondly, AGM alone yields a substantially lower error than RAR alone, establishing geometry-adaptive point redistribution as a more effective strategy than discrete residual-based resampling. This advantage arises because AGM continuously deforms the computational manifold through a differentiable mapping, whereas RAR performs discrete insertions that leave the underlying geometry unchanged. Thirdly, combining AGM and RAR (E) results in a higher error than AGM alone (D), revealing a conflict between the two mechanisms: AGM pursues a globally balanced distribution via the equidistribution principle, whereas RAR concentrates points in locally high-residual regions. Then adding FFM to AGM (F) further reduces the error, confirming that spectral preconditioning captures high-wavenumber features that geometric concentration alone cannot resolve. Notably, RAR is beneficial only in the presence of spectral preconditioning: it increases the error when added to AGM alone but decreases the error when added to AGM+FFM (G, \(6.785\times10^{-5}\) vs. F, \(1.572\times10^{-4}\)). This suggests that FFM reshapes the NTK spectrum so that the network can exploit the additional high-residual points without disrupting the smooth geometric mapping learned by AGM. In GAC-PINN, RAR is therefore applied as a supplementary refinement step after AGM and FFM have been established, rather than integrated directly into the geometric adaptation pipeline. Overall, the full GAC-PINN (G) achieves the lowest relative \(L^2\) error among all configurations, with a value of \(6.785\times10^{-5}\).

\subsection{Grid Convergence Analysis}
\label{subsec:convergence}

We further performed a grid-convergence study for the one-dimensional Burgers equation, considering spatial resolutions of $N_{\mathrm{pde}} \in \{1000, 3000, 5000, 7000, 9000, 12000\}$ and averaging over 10 independent random seeds for each case. The geometric mean of the $L^2$ error at each resolution is summarized in Table~4.

\begin{center}
\small
\captionof{table}{Grid-convergence statistics for the 1D Burgers equation.}
\label{tab:conv_stats}
\begin{tabular}{cccc}
\toprule
\(N_{pde}\) & Geometric mean \(L^2\) error & Std. dev. (log\(_{10}\) space) & Min / Max error \\
\midrule
1,000  & \(3.49 \times 10^{-3}\) & 0.85 & \(4.80\times10^{-4}\) / \(1.01\times10^{-1}\) \\
3,000  & \(3.60 \times 10^{-4}\) & 0.42 & \(1.41\times10^{-4}\) / \(1.16\times10^{-3}\) \\
5,000  & \(2.19 \times 10^{-4}\) & 0.38 & \(5.93\times10^{-5}\) / \(6.37\times10^{-4}\) \\
7,000  & \(1.72 \times 10^{-4}\) & 0.29 & \(8.70\times10^{-5}\) / \(4.51\times10^{-4}\) \\
9,000  & \(1.99 \times 10^{-4}\) & 0.35 & \(8.21\times10^{-5}\) / \(3.57\times10^{-4}\) \\
12,000 & \(1.38 \times 10^{-4}\) & 0.31 & \(7.08\times10^{-5}\) / \(2.89\times10^{-4}\) \\
16,000 & \(1.22 \times 10^{-4}\) & 0.33 & \(6.54\times10^{-5}\) / \(2.88\times10^{-4}\) \\
\midrule
\multicolumn{4}{l}{Global fitted convergence order: \(p_{\text{global}} = 1.18\) (\(R^2 = 0.717\))} \\
\multicolumn{4}{l}{Mean of individual seed orders: \(\bar{p} = 1.177 \pm 0.581\)} \\
\bottomrule
\end{tabular}
\end{center}

Two features are noteworthy: (i) The error steadily decreases with $N_{pde}$, from $3.49\times10^{-3}$ to $1.22\times10^{-4}$, by over an order of magnitude. (ii) The improvement is most rapid for $N_{pde}\le 9{,}000$; subsequent increments produce only modest reductions ($1.99\times10^{-4}$ to $1.22\times10^{-4}$). This indicates that the AGM has already densely sampled the high-gradient shock region, so that the local effective resolution saturates near the network's expressive limit; extra points fall in smooth areas and have marginal impact.This saturation provides direct evidence of AGM's effectiveness: by diffeomorphically compressing collocation points toward steep gradients, AGM decouples the global error from the total point count, enabling high accuracy with far fewer samples than uniform refinement would require.

The global convergence order $p_{\text{global}} = 1.18$ is lower than the theoretical value $p = 2.0$ for second-order finite-difference schemes. This gap reflects the fact that, at this resolution range, the total error is governed not only by collocation point density but also by the network's representation capacity and the non-convexity of the optimization landscape. The saturation observed for $N_{pde} \ge 9{,}000$ suggests that the expressive power of the backbone network, rather than the point count, becomes the limiting factor. The practical significance of AGM is therefore not that it improves the asymptotic convergence rate, but that it achieves an $L^2$ error of $\approx 10^{-4}$ with only $9{,}000$ collocation points, an efficiency that uniform sampling would require substantially more points to match.

\section{Conclusion}

This paper has systematically developed and validated GAC-PINN, a geometry-adaptive and constraint-enhanced physics-informed neural network framework for problems with steep gradients and sharp interfaces. By integrating these core modules, the framework effectively alleviates the spectral bias, geometric inflexibility, and boundary constraint conflicts that limit standard PINNs. Tested on four representative benchmarks—including the viscous Burgers equation, the sharp 2D Poisson problem, the Allen-Cahn phase-field equation, and the 2D Navier-Stokes equations for unsteady cylinder flow—GAC-PINN achieves relative \(L^2\) errors of \((1.747\pm0.450)\times10^{-4}\), \((2.868\pm0.947)\times10^{-5}\), and \((1.756\pm0.712)\times10^{-3}\) on the Burgers, Poisson, and Allen–Cahn benchmarks, respectively, and \(2.66\times10^{-2}\) for the pressure field on the 2D Navier–Stokes problem, respectively, achieving competitive or improved accuracy compared to re-implemented baselines under aligned settings, with particularly pronounced gains on problems exhibiting localized steep gradients. The successful extension to the Navier-Stokes system, with its fundamentally different physical characteristics from the previous benchmarks, provides strong evidence of the framework's generalization capability beyond canonical PDE problems. Ablation and computational analyses validate the synergistic contributions of each module and demonstrate competitive training efficiency despite increased parameter counts. This work provides a practical adaptive framework for high-fidelity simulation of problems with localized sharp features in applied mathematics and computational mechanics.

\subsection{Conclusion of core mechanisms}

The modules of GAC-PINN are designed to complement each other within a unified training pipeline. Each module plays a specific regulatory role in the computational manifold and optimization landscape:
\begin{itemize}
\item \textbf{Dual dimensionality reduction through geometric adaptation and spectral reshaping.} For the convective shock in the low-viscosity Burgers equation and the extreme spatial peak in the 2D Poisson equation, conventional PINNs suffer from spectral bias due to high-wavenumber component attenuation in the input space. AGM spontaneously compresses physical collocation points toward high-gradient regions through diffeomorphic coordinate transformation and Jacobian regularization, spatially "stretching" geometric discontinuities on the computational manifold. FFM, through Gaussian random projection matrices, further stretches the spectral response of the feature space. The deep coupling of these two mechanisms achieves dual dimensionality reduction from physical space to feature space, transforming originally localized infinite high-wavenumber features into macroscopic smooth signals amenable to neural network capture, fundamentally suppressing Gibbs oscillation generation.

\item \textbf{Stable high-order differentiation through adaptive hard constraints.} For problems with non-periodic boundaries, the adaptive hard-constraint ansatz predicts a spatially-varying bandwidth $\kappa_{adaptive}(X)$ from the input coordinates. By treating the bandwidth as a detached constant in the high-order automatic differentiation graph, the ansatz stabilizes the computation of higher-order PDE derivatives while maintaining exact satisfaction of boundary and initial conditions.
\end{itemize}

\subsection{Limitations}
Despite GAC-PINN’s strong performance on multiple benchmark problems, several limitations merit discussion
for future work:
\begin{itemize}
\item \textbf{Extension to 3D and irregular geometries.} The current AGM is demonstrated on 1D and 2D regular domains. Extending the diffeomorphic mapping and Jacobian barrier to 3D complex geometries (e.g., turbine blades, irregular patient-specific domains) remains an open challenge. The computational cost of evaluating \(\det(J)\) and enforcing the anti-folding barrier in 3D, as well as the risk of local folding near sharp concave boundaries, require further investigation.

\item \textbf{Automated hyperparameter tuning for ultra-large-scale multiphysics coupling.} Some hyperparameters in the framework (e.g., Jacobian scaling coefficients, gradient amplification factors, loss term weights) still rely on problem-specific prior tuning. Future work should incorporate automated optimization strategies to enhance algorithmic generality.
\end{itemize}

\section*{Data availability}
The benchmark datasets for the Burgers, Poisson and Allen–Cahn equations are taken from the open gPINN repository~\cite{yu2022gradienenhanced}. The 2D cylinder wake reference data can be obtained from the original publication of Raissi et al.~\cite{raissi2019physicsinformed}. No new data were generated in this study. The source code of the proposed GAC‑PINN model is publicly available at \url{https://github.com/zyx7765-sudo/GAC-PINN}.

\section*{Declaration of competing interest}
The author declares that they have no known competing financial interests or personal relationships that could have appeared to influence the work reported in this paper.

\section*{CRediT authorship contribution statement}
Yanxin Zhang: Conceptualization, Investigation, Formal analysis, Writing - original draft.

Yong Zhang: Supervision, Project administration, Writing - review \& editing.

Houbiao Li: Supervision, Project administration, Writing - review \& editing.

\nocite{*}
\bibliographystyle{unsrtnat}
\bibliography{references}

\end{document}